\documentclass[12pt]{article}
\usepackage{fullpage}
\usepackage{amssymb}
\usepackage{amsthm}
\usepackage[all,2cell]{xy}
\newcommand{\ncom}{\newcommand}
\ncom{\beqn}{\begin{eqnarray*}}
\ncom{\eeqn}{\end{eqnarray*}}

\newtheorem{theorem}{Theorem}[section]
\newtheorem{corollary}[theorem]{Corollary}

\newtheorem{lemma}[theorem]{Lemma}

\newtheorem{proposition}[theorem]{Proposition}
\newtheorem{remark}[theorem]{Remark}

\begin{document}
\title{Torus quotients of homogeneous spaces of the general linear group 
and the standard representation of certain symmetric groups.}
\author{S.S.Kannan, Pranab Sardar  \\  
\\ Chennai Mathematical Institute, Plot H1, SIPCOT IT Park,\\ Padur 
Post Office, Siruseri, Tamilnadu - 603103, India.\\
kannan@cmi.ac.in, pranab@cmi.ac.in} 

\maketitle
\date{14.08.2007}

\begin{abstract} We give a stratification of the GIT quotient of the 
Grassmannian $G_{2,n}$ modulo the normaliser of a maximal torus of 
$SL_{n}(k)$ with respect to the ample generator of the Picard group 
of $G_{2,n}$. We also prove that the flag variety $GL_{n}(k)/B_{n}$ can be 
obtained as a GIT quotient of $GL_{n+1}(k)/B_{n+1}$ modulo a maximal torus of 
$SL_{n+1}(k)$ for a suitable choice of an ample line bundle on 
$GL_{n+1}(k)/B_{n+1}$.    
\end{abstract}
\hspace*{2cm}Keywords: GIT quotient, line bundle, simple reflection.

\begin{center}
{\bf Introduction}
\end{center}

 Let $k$ be an algebraically closed field. Consider the action of a maximal 
torus $T$ of $SL_{n}(k)$ on 
the Grassmannian $G_{r,n}$ of $r$- dimensional vector subspaces of an $n$- 
dimensional vector space over $k$. Let $N$ denote the normaliser of T in 
$SL_n(k)$. Let $\mathcal{L}_r$ 
denote the ample generator of the Picard group of $G_{r,n}$. Let $W=N/T$ 
denote the Weyl group of $SL_{n}(k)$ with respect to $T$.

In [5], it is shown that the semi-stable points of $G_{r,n}$ with respect to 
the $T$-linearised line bundle $\mathcal L_r$ is same as the stable points if 
and only if $r$ and $n$ are co-prime.

In this paper, we describe all the semi-stable points of $G_{r,n}$ with respect
to $\mathcal L_r$. In this connection, we prove the following result:\\
First, we introduce some notation needed for the statement of the theorem.

Let $\mathfrak{h}_{j}$ be a Cartan subalgebra of 
$\mathfrak{sl}_{j+1}$, $\mathbb P(\mathfrak{h}_j)$ be the projective space and 
$R_{j}\subseteq \mathfrak{h}_{j}^{*}$ be the root system. Let $V_j$ be the open 
subset of  $\mathbb P(\mathfrak{h}_{j})$ defined by 
\[V_j:=\{x\in \mathbb P(\mathfrak{h}_{j}):\alpha(x)\neq 0, \forall \alpha 
\in R_{j}\}.\]  Here, the Weyl group of $\mathfrak{sl}_{j+1}$ is $S_{j+1}$, and 
$\mathfrak{h}_{j}$ is the standard representation of $S_{j+1}$. 

With this notation, taking $m= \lceil{\frac{n-1}{2}}\rceil$ (for this 
notation, see lemma  1.6) and 
$t=[\frac {n-1}{2}]$ we have \\
{\bf Theorem:}
  ${_{_{N}\backslash \backslash}}G_{2,n}^{ss}(\mathcal{L}_2)$  has a 
stratification $\bigcup_{i=0}^{t}C_i$ where 
$C_0= {_{S_{m+1}}\backslash}{\mathbb P(\mathfrak{h}_m)}$, and $C_i=
{_{S_{i+m+1}\backslash}{V_{i+m}}}$. 

On the other hand, the GIT quotient of $GL_{n+1}(k)/B_{n+1}$ modulo a maximal 
torus of $SL_{n+1}(k)$ for any ample line bundle on $GL_{n+1}(k)/B_{n+1}$ 
and $GL_n(k)/B_n$ are both birational varieties. So, it is a natural question 
to ask whether the flag variety $GL_{n}(k)/B_{n}$ can be 
obtained as a GIT quotient of $GL_{n+1}(k)/B_{n+1}$ modulo a maximal torus of 
$SL_{n+1}(k)$ for a suitable choice of an ample line bundle on 
$GL_{n+1}(k)/B_{n+1}$. We give an affirmative answer to this question.
For a more precise statement, see theorem 5.2. In this connection, we also 
prove that the action of the Weyl group $S_{n+1}$ on the quotient is given by 
the standard representation. For a more precise statement, see corollary 5.4. 

Section 1 consists of preliminary notation and some combinatorial lemmas about 
minuscule weights.

In section 2, we describe all Schubert cells in $G_{r,n}$ admitting semi-stable 
points.

In section 3, we describe the action of the Weyl group $W$ on  
${_{_{T}\backslash \backslash}}G_{r,n}^{ss}(\mathcal{L}_r)$.

In section 4, we describe a stratification of  
${_{_{N}\backslash \backslash}}G_{2,n}^{ss}(\mathcal{L}_2)$.

In section 5, we obtain  $GL_{n}(k)/B_{n}$ as a GIT quotient of 
$GL_{n+1}(k)/B_{n+1}$ modulo a maximal torus of $SL_{n+1}(k)$
for a suitable line bundle on $GL_{n+1}(k)/B_{n+1}$.

\section{Preliminary notation and some combinatorial Lemmas}
This section consists of preliminary notation and some combinatorial lemmas 
about minuscule weights.
Let $G$ be a reductive Chevalley group over an algebraically closed field $k$. 
Let $T$ be a 
maximal torus of the commutator subgroup $[G, G]$, $B$ a Borel subgroup of 
$G$ containing $T$ and $U$ be the unipotent radical of $B$. 
Let $N$ be the normaliser of $T$ in $[G,G]$. Let $W=N/T$ be Weyl group of 
$[G, G]$ with respect to $T$ and $R$ denote the set of roots with respect to 
$T$, $R^{+}$ positive roots with respect to $B$. Let $U_{\alpha}$ denote 
the one dimentional $T$-stable subgroup of $G$ corresponding to the root 
$\alpha$ and let $S=\{\alpha_1, \cdots \alpha_l\}\in R^{+}$ denote the set of 
simple roots.  
For a subset $I\subseteq S$ denote 
$W^I =\{w\in W| w(\alpha)>0, \, \alpha\in I\}$.
Let $X(T)$ (resp. $Y(T)$) denote the set of characters of $T$ (resp. 
one parameter subgroups of $T$). Let $E_1:= X(T)\otimes \mathbb R$, 
$E_2=Y(T)\otimes \mathbb R$. Let $\langle . ,  .\rangle : E_1\times E_2 
\longrightarrow \mathbb R$ be the canonical non-degenerate bilinear form. 
Choose $\lambda_j$'s in $E_2$ such that 
$\langle \alpha_i, \lambda_j \rangle = \delta_{ij}$ for all $i$. Let 
$\overline{C(B)}:=\mathbb R _{\geq 0}$ - span of the $\lambda_i$'s . 
Let $\check{\alpha} \in Y(T)$ be as in page-19 of [1]. We also have 
$s_{\alpha}(\chi)=\chi-\langle \chi, \check{\alpha}\rangle \alpha$ for all 
$\alpha \in R$ and $\chi \in E_1$. Set $s_i=s_{\alpha_i} \,\ \forall 
\,\ i=1,2 \cdots l$. Let $\{\omega_i: i=1,2\cdots l\} \subset E_1$ be the 
fundamental weights; i.e. $\langle \omega_i, \check{\alpha_j} \rangle = 
\delta_{ij}$ for all $i,j = 1,2 \cdots l$.  

We now prove some elementary lemmas about minuscule weights. For 
notation, we refer to [7].   

\begin{lemma}
Let $I$ be any nonempty subset of $S$, and let $\mu$ be a weight of the 
form $\sum_{\alpha_i \in I}m_i\alpha_i -\, \sum_{\alpha_i\not\in I}m_i\alpha_i$,
where $m_i\in \mathbb Q$ for all $i$, $1\leq i\leq l$; $m_i>0$ for all 
$\alpha_i \in I$ and $m_i\geq 0$ for all $\alpha_i \in S\setminus I$.
Then there is an $\alpha \in I$ such that $s_{\alpha}(\mu)<\mu$.
\end{lemma}
\begin{proof}
Since 
$s_{\alpha}(\mu)=\,\mu-\langle \mu, \check{\alpha}\rangle\alpha$, we need 
to find an $\alpha \in I$ such that $\langle\mu, \check{\alpha}\rangle >0$. 
This follows because the Cartan matrix $(\langle \alpha_i, 
\check{\alpha_j}\rangle)_{i,j}$
is positive definite, so we can find an $\alpha \in I$ such that 
$\langle\sum_{\alpha_i \in I}m_i\alpha_i, \check{\alpha}\rangle >0$. Now we 
know that for any  
$\alpha_i, \alpha_j \in S$, $i\neq j$, $\langle\alpha_i,\check{\alpha_j}\
rangle\leq 0$. Hence,
$\langle \sum_{\alpha_i\not\in I}m_i\alpha_i,\check{\alpha}\rangle \leq 0$
for this $\alpha \in I$. Thus $\langle\mu, \check{\alpha}\rangle>0$. 
This proves the lemma. 
\end{proof}

\begin{lemma}
Let $\lambda$ be any dominant weight and let 
$I=\{\alpha\in S: \langle\lambda,\check{\alpha}\rangle =0\}$.
Let $w_1, w_2 \in W^I$ be such that $w_1(\lambda)=\, w_2(\lambda)$. 
Then $w_1=w_2$.
\end{lemma}
\begin{proof}
See [1] and [2].
\end{proof}

\emph{In the rest of this section, $\omega$ will denote a minuscule weight
and $I:= \{\alpha\in S:\langle \omega, \check{\alpha}\rangle =0\}$}
\begin{lemma}
Let $\alpha \in S$ and $\tau \in W$ such that 
$l(s_{\alpha} \tau)=l(\tau)+1$ and $ s_{\alpha} \tau \in W^I$, then
$\tau \in W^I$; $s_{\alpha} \tau(\omega)=\tau(\omega)-\alpha$.
\end{lemma}
\begin{proof}
The proof of the first part of the lemma is clear. Now 
$s_{\alpha} \tau(\omega)= \tau(\omega)-\langle\tau(\omega),\check{\alpha}
\rangle\alpha$.
Since the pairing $\langle.,.\rangle$ is $W$-invariant, 
$\langle\tau(\omega),\check{\alpha}\rangle\,=\,\langle\omega,
\check{\tau^{-1}\alpha}\rangle$. Again since $l(s_{\alpha}\tau)=l(\tau)+1$, 
we have $\tau^{-1}\alpha>0$. 
Let $\check{\tau^{-1}\alpha}=\sum_{i=1}^{l}m_i\check{\alpha_i}$, 
$m_i\in \mathbb Z_{\geq 0}$.
Now, if $\langle\omega,\check{\tau^{-1}\alpha}\rangle=0$, then 
$m_i>0\Rightarrow \langle\omega,\check{\tau^{-1}\alpha_i}\rangle=0$ for 
$1\leq i\leq l$. This gives a contradiction, since $s_{\alpha}\tau \in W^I$
and $s_{\alpha}\tau(\tau^{-1}\alpha)=s_{\alpha}(\alpha)<0$.
Thus, $\langle\omega,\check{\tau^{-1}\alpha}\rangle =\, 1$. Hence the lemma is
proved.
\end{proof}

\begin{corollary}
$1$. For any $w\in W^I$, the number of times that $s_i$, $1\leq i\leq n-1$
appears in a reduced expression of 
$w\,=$ (coefficient of $\alpha_i$ in $\omega)\,-$ 
(coefficient of $\alpha_i$ in $w(\omega)$) and hence it is independent
of the reduced expression of $w$.\\
$2$. Let $w\in W^I$ and let $w=s_{i_1}.s_{i_2}\ldots s_{i_k}\in W^I$ be a
reduced expression. Then $w(\omega)=\omega-\sum_{j=1}^{k}\alpha_{i_j}$. 
and $l(w)=ht(\omega -w(\omega))$.\\
\end{corollary}
\begin{proof}
Follows from Lemma 1.3. 
\end{proof}

\begin{lemma}
Let $w=s_{i_1}s_{i_2}\ldots s_{i_k}\in W$ such that 
$ht(\omega- s_{i_1}s_{i_2}\ldots s_{i_k} (\omega))=k$ then
$w\in W^I$ and $l(w)=k$.
\end{lemma}
\begin{proof}
This follows from the corollary 1.4.
\end{proof}

\begin{lemma}
Let $\omega=\sum_{i=1}^{l}m_i\alpha_i$, $m_i\in \mathbb Q_{\geq 0}$
be a minuscule weight. 
Let $I=\{\alpha\in S: \langle\omega,\check {\alpha}\rangle=0\}$.
Then, there exist a unique $w\in W^I$ such that 
$w(\omega)=\sum_{i=1}^{l}(m_i-\lceil m_i \rceil)\alpha_i$ where for any
real number $x$, 
\[
\lceil x \rceil := \,\cases{ x & if x is an integer \cr
                       [x]+1 & otherwise}
\]  
\end{lemma}
\begin{proof}
Using lemma $1.1$ and the fact that $\omega$ is minuscule we 
can find a sequence  $s_{i_k},s_{i_{k-1}}, \cdots \\,s_{i_1}$ 
of simple reflections in $W$ such that 
for each $j$, $2\leq j\leq k+1$, coefficient of $\alpha_{i_j}$ in 
$s_{i_{j-1}}.s_{i_{j-2}}\ldots s_{i_1}(\omega_r)$ is positive 
 and $(s_{i_k}.s_{i_{k-1}}\ldots 
s_{i_1}(\omega_r))=\,\omega_r- \sum_{j=1}^{k}\alpha_{i_j}$ for each 
$j$, $1\leq j\leq k$.
The existence part of the lemma follows from here. The uniqueness
follows from lemma $1.2$.
\end{proof}

\begin{lemma}
Let $\omega=\sum_{i=1}^{l}m_i\alpha_i$, $m_i\in \mathbb Q_{\geq 0}$
be a minuscule weight. Let $I=\{\alpha\in S: \langle\omega,\check{\alpha}
\rangle=0\}$.
Then, there exist a unique $\tau\in W^I$ such that 
$\tau(\omega)=\sum_{i=1}^{l}(m_i- [m_i])\alpha_i$.  
\end{lemma}
\begin{proof}
Proof is similar to that of lemma 1.6.
\end{proof}

Now onwards, we say that for two elements $w$ and $\tau$ in $W$, $w\leq \tau$
if $l(\tau)=l(w)+l(\tau w^{-1})$.

\begin{lemma}
Let $\omega$ and $I$ be as in the lemma $1.6$ and $\tau,\sigma\in W^I$.
Then $\tau(\omega)\leq \sigma(\omega) \Leftrightarrow \sigma \leq \tau$.
\end{lemma}

\begin{proof}
The proof is by induction on
$ht(\sigma(\omega)-\tau(\omega))$ which is a non-negative integer.\\
$\underline{ht(w(\sigma\omega)-\tau(\omega))=1}$: 
This means $\sigma(\omega)= \tau(\omega) + \alpha$ for some
$\alpha \in S$. Applying $s_{\alpha}$ on both the sides of this
equation, we have,
\[ 
\begin{array}{l}
s_{\alpha}\sigma(\omega)=\,-\alpha +s_{\alpha}\tau(\omega)\\
\Longrightarrow \tau(\omega)-\langle \omega,\check{\sigma^{-1}\alpha}\rangle 
\alpha\,=
-2\alpha +\tau(\omega)-\langle \omega, \check{\tau^{-1}\alpha} \rangle \alpha \\
\Longrightarrow \langle \omega,\check{\sigma^{-1}\alpha} \rangle=\,2+\langle \omega, 
\check{\tau^{-1}\alpha} \rangle
\end{array}
\]

Since  $\omega$ is minuscule, we get 
$\langle \omega,\check{\sigma^{-1}\alpha} \rangle= 1$ and 
$\langle \omega, \check{\tau^{-1}\alpha} \rangle=-1$. This implies, by the
lemma $1.5$, that $l(s_{\alpha} \sigma)= l(w)+1$ and $s_{\alpha} w\in W^I$.
Now, we have $s_{\alpha}\sigma(\omega)=\tau(\omega)$. Hence, by
lemma $1.2$, we get $\tau=s_{\alpha}\sigma$ with $l(\tau)=l(\sigma)+1$. 
Thus the result follows in this case. \\
Let us assume that the result is true for 
$ht(\sigma(\omega)-\tau(\omega))\leq m-1$. \\
$\underline{ht(\sigma(\omega)-\tau(\omega))=m}$: 
Let $\sigma(\omega)-\tau(\omega)=\sum_{\alpha_i\in J}m_i\alpha_i$ where
$J\subseteq S$ and $m_i$'s are positive integers. Since $\langle 
\sum_{\alpha_i\in J}m_i\alpha_i, \sum_{\alpha_i\in J}m_i\check{\alpha_i}\rangle
\geq 0$ there exist an $\alpha_j \in J$ such that 
$\langle \sigma(\omega)-\tau(\omega),\check{\alpha_j}\rangle >0$. Hence either
$\langle \sigma(\omega),\check{\alpha_j}\rangle >0 $ or 
$\langle \tau(\omega), \check{\alpha_j}\rangle <0$. \\
\underline{Case $I$}: Let us assume $\langle \sigma(\omega),
\check{\alpha_j}\rangle >0$ . 
Then $l(s_{\alpha_j}\sigma)=l(\sigma)+1$ and $s_{\alpha_j} \sigma \in W^I$. 
Now $ht(s_{\alpha_j}\sigma (\omega)-\tau(\omega))=m-1$. Hence, by induction
$\tau=\phi_1s_{\alpha_j}\sigma$ with $l(\tau)=l(\phi_1)+l(s_{\alpha_j}\sigma)$.
Thus taking $\phi=\phi_1.s_{\alpha_j}$ we are done in this case.\\
\underline{Case $II$}: Let us assume $\langle \tau(\omega), 
\check{\alpha_j}\rangle <0$. Then $l(s_{\alpha_j}\tau)=l(\tau)-1$ and 
$s_{\alpha_j}\tau \in W^I$. Since 
$\sigma(\omega)-s_{\alpha_j}\tau(\omega)=m-1$ by induction 
$s_{\alpha_j}\tau=\phi_2\sigma$ with $l(s_{\alpha_j}\tau)=l(\phi_2)+l(\sigma)$.
Thus taking $\phi=s_{\alpha_j}\phi_2$ we are done in this case also.
This completes the proof. 
\end{proof}

\begin{corollary}
Let $\omega$, $w$ and $I$ be as in lemma 1.6. Let $\sigma \in W^I$ be such that
$\sigma(n\omega)\leq 0$ for some positive integer. Then, we have 
$w \leq \sigma$.
\end{corollary}
\begin{proof}
The proof follows from lemma $1.6$, $1.8$ and the fact that $\omega$
is minuscule.
\end{proof}

\begin{corollary}
Let $\omega$, $w$ and $I$ be as in lemma 1.6. Let $\sigma \in W^I$ be such that
$\sigma(n\omega)\geq 0$ for some positive integer. Then, we have 
$\sigma \leq w$
\end{corollary}
\begin{proof}
The proof follows from lemma $1.7$, $1.8$ and the fact that $\omega$
is minuscule.
\end{proof}

\section{Description of Schubert varieties  in the Grassmannian having semi-stable points}
In this section, we have the following notation. Let $G= GL_{n}(k)$ with 
characteristic of $k$ is either zero or bigger than $n$. Let $r\in \{2, 
\cdots n-2\}$.
Consider the action of a maximal torus $T$ of $SL_{n}(k)$ on the Grassmannian  
$G_{r,n}$. Let $B$ be a Borel subgroup  of $G$ containing $T$. Let 
$S=\{\alpha_1, \cdots \alpha_{n-1}\}$ be the set of simple roots with respect 
to $B$ arranged in the ordering of the vertices in the Dynkin diagram of type 
$A_{n-1}$. Let 
$I_r=S\setminus \{\alpha_{r}\}$. We first note that $G_{r,n}$ is the 
homogeneous space $GL_{n}(k)/P_r$ where $P_r=BW_{I_r}B$ is the maximal 
parabolic subgroup of $GL_{n}(k)$ containing $B$ associated to the simple root 
$\alpha_{r}$. Let $\omega_{r}$ be the fundamental weight associated to the 
simple root $\alpha_{r}$ and let $\mathcal{L}_r$ denote the line bundle on 
$GL_{n}(k)/P_r$ corresponding to $\omega_{r}$. We describe all Schubert cells 
in $GL_{n}(k)/P_r$ admitting semi-stable points for the above mentioned action 
of $T$ with respect to the line bundle $\mathcal{L}_r$.

Some of the elementary facts about the combinatorics of $W^{I_{r}}$ that is 
being used in this section can be found in [7]. For the 
convenience of the reader, we prove them here.

\begin{lemma}
Let $w\in W^I, w\neq id$. Then there exists an $i\in \mathbb N$, $i\leq r$
and a sequence of positive integers $\{a_j\}$, $j=1,2,\ldots, r$ such that
the following holds.\\
$(a)$ $a_j\geq j$ for all $j$, $i\leq j\leq r$\\
$(b)$ 
$w=(s_{a_i}.s_{a_i-1}\ldots s_i)(s_{a_{i+1}}.s_{a_{i+1}-1}\ldots s_{i+1})
\ldots(s_{a_r}.s_{a_r-1}\ldots s_r)$ with $l(w)=\sum_{j=i}^{r}(a_j-j+1)$
\end{lemma}

\begin{proof}
Let $i$ be the least positive integer such that $s_{\alpha_i}\leq w$.
The rest of the proof follows from braid relations in $W$.
\end{proof}

\begin{lemma}
Let $w,\tau\in W^I$. Write 
$w=(s_{a_i}.s_{a_i-1}\ldots s_i)(s_{a_{i+1}}.s_{a_{i+1}-1}\ldots s_{i+1})
\ldots(s_{a_r}.s_{a_r-1}\ldots s_r)$
and 
$\tau=(s_{b_k}.s_{b_k-1}\ldots s_k)(s_{b_{k+1}}.s_{b_{k+1}-1}\ldots s_{k+1})
\ldots(s_{b_r}.s_{b_r-1}\ldots s_r)$ be
as in the lemma $2.1$. Then $w\leq \tau \Leftrightarrow k\leq i$ and 
$b_j\geq a_j$ for all $j$, $i\leq j\leq r$.
\end{lemma}

\begin{proof} The proof follows from lemma $1.8$ and the fact 
that $w(\omega_r)\geq \tau(\omega_r) \Leftrightarrow k\leq i$ and 
$b_j\geq a_j$ for all $j$, $i\leq j\leq r$.
\end{proof}

Now, write $n=qr+t$ with $1\leq t \leq r$ and  let $\tau_{r}\in W^{I_{r}}$ be 
the unique element as in lemma 1.6 for 
the case when $\omega=\omega_{r}$. Then, $\tau_{r}$ must be of the form 
$\tau_{r}=(s_{a_{1} }\cdots s_1)\cdots (s_{a_{r}}\cdots s_{r})$ where 

$$a_i= \left \{ \begin{array}{l} i(q+1) ~ ~ if  ~ ~ i \leq t-1.\\ 

iq+(t-1) ~ ~ if ~ ~ t\leq i\leq r\\ 

\end{array} \right .$$  

Let $\tau^{n-r}\in W^{I_{n-r}}$ be the unique element as in lemma 1.7 for the case
$\omega= \omega_r$. Then, we have $\tau_{r}=\tau^{n-r}w_{0}^{I_{r}}$ and 
$l(w_{0}^{I_{r}})=l(\tau_{r})+ l(\tau^{n-r})$. 

Let $w\in W^I$ be such that $w(n\omega_r)\leq 0$. 

Then, we have 

\begin{lemma} $\tau_{r}\leq w$ and $w \tau_{r}^{-1}\leq (\tau^{n-r})^{-1}$.
\end{lemma} 
\begin{proof}
Proof follows from corollary 1.8 and corollary 1.9.
\end{proof}

For any such $w$, we describe the set $R^{+}(w^{-1})$.

\begin{lemma}
$R^{+}(w^{-1})$ consists of roots of the form 
$\alpha_j+\alpha_{j+1}+\cdots +\alpha_{a_i}$ for $1\leq i\leq r$ where 
$j\neq a_k+1$ for any $k<i$.
\end{lemma}

\begin{proof}
We have 
$w^{-1}=(s_r\ldots s_{a_r})\ldots (s_2\ldots s_{a_2}).(s_1\ldots s_{a_1})$,
which is a reduced expression. Thus the elements of $R^{+}(w^{-1})$ are
\[\beta_{i,j-i+1}=(s_{a_1}\ldots s_1).(s_{a_2}\ldots s_2)\ldots (s_{a_i}
\ldots s_{j+1}.\hat{s_j}.\hat{s_{j-1}}\ldots \hat{s_i})(\alpha_j) \]
where $i\leq j\leq a_i$, $1\leq i\leq r$, $\hat{}$ denotes omission 
of the symbols. We have, 
\[(s_{a_i}\ldots s_{j+1}.\hat{s_j}.\hat{s_{j-1}}\ldots \hat{s_i})(\alpha_j)
=\alpha_{j}+\alpha_{j+1}+\cdots +\alpha_{a_i}\]
Since, $a_1<a_2<\cdots< a_r$, each $\beta_{i,j}$ is of the form
\[ \alpha_{j}+\alpha_{j+1}+\cdots +\alpha_{a_i}. \] Now $j\neq a_{k}+1$ for any
$k < i$ follows from the fact that $l(w)$ is the same as the cardinality of 
$R^{+}(w^{-1})$.
\end{proof}

\begin{remark}
From the lemma it follows that the elements of $R^{+}(w^{-1})$ can be written
in an array as follows:
\[
\begin{array}{lllllllllllllll}
\beta_{1,1}& \beta_{1,2}&\cdots & \beta_{1,a_1}&&&&&&&&&&& \\
\beta{2,1} &\beta_{2,2} & \cdots &\beta_{2,a_1}& \beta_{2,a_1+1}&\beta_{2,a_1+2}& \cdots & \beta_{2,a_2-1}&&&&&&& \\
\beta{3,1} &\beta_{3,2} & \cdots &\beta_{3,a_1}& \beta_{3,a_1+1}&\beta_{3,a_1+2}& \cdots & \beta_{3,a_2-1} & \beta_{3,a_2}&\cdots & \beta_{3,a_3-2}&& &&\\
\vdots & \vdots & &\vdots &\vdots &\vdots && \vdots & \vdots & & \vdots && &&\\
\beta{r,1} &\beta_{r,2} & \cdots &\beta_{r,a_1}& \beta_{r,a_1+1}&\beta_{r,a_1+2}& \cdots & \beta_{r,a_2-1} & \beta_{r,a_2}&\cdots & \beta_{r,a_3-2}& \cdots & \beta_{r,a_r-r+1} && \\

\end{array}
\]
where the array has $r$ rows, and the length of the $i$-th row is $a_i-(i-1)$. 
Note that $\beta_{1,a_1}=\alpha_{a_1}$, and for $2\leq i\leq r$, 
$\beta_{i,a_i-i+1}=\alpha_{a_i}$, only if $a_i\geq a_{i-1}+2$. 
In this case, for all $j$, $i\leq j\leq r$, 
$\beta_{j,a_{i-1}-i+2}=\beta_{i-1,a_{i-1}-i+2}+\alpha_{a_{i-1}+1}+\alpha_{a_{i-1}+2}+
\cdots +\alpha_{a_j}$ and 
$\beta_{j,a_{i-1}-i+3}=\alpha_{a_{i-1}+2}+\alpha_{a_{i-1}+3}+\cdots +\alpha_{a_j}$.
If $a_i=a_{i-1}+1$, then $a_i-i+1=a_{i-1}-(i-1)+1$, therefore, the $(i-1)$-th 
and $i$-th rows have same length. In this case for all $j$, $i\leq j\leq r$, 
$\beta_{j,a_i-i+1}=\beta_{i-1,a_i-i+1}+ \alpha_{a_{i-1}+1}+\alpha_{a_{i-1}+2}+\cdots 
+ \alpha_{a_j}$. 
\end{remark}

For any $w\in W^{I}$, let $X(w):=\overline{BwP_r/P_r}$ denote the Schubert 
variety in $GL_{n}(k)/P_r$.

We recall $BwP_{r}/P_{r}=U_{w}wP_{r}$, where $U_{w}$ is the product 
$\prod_{\alpha\in R^{+}(w^{-1})}U_{\alpha}$ of the root groups $U_{\alpha}$, 
and we describe below the ordering of roots in which the product is taken.

Consider the open set 
\[ V:=\{\prod_{\beta_{ij}\in R^{+}(w^{-1})}u_{\beta_{ij}}(x_{\beta_{ij}}).w.P_r:
x_{\beta_{ij}}\neq 0,\,\forall \beta_{ij}\in R^{+}(w^{-1})\} \]
of $X(w)$ in $GL_n/P_r$ where the order in which the
product is taken is as follows: Put a partial order on $R^{+}(w^{-1})$ by
declaring $\beta_{ij}\leq \beta_{kl}$ if either $i=k$ and $j\geq l$ or if
$i<k$. Now we take the product so that whenever $\beta_{ij}\leq \beta_{kl}$,
$u_{\beta_{ij}}(x_{\beta_{ij}})$ appears on the right hand side of 
$u_{\beta_{kl}}(x_{\beta_{kl}})$. Note that $u_{\beta_{ij}}(x_{\beta_{ij}})$'s commute
with each other, since $\beta_{i_1,j_1}, \beta_{i_2,j_2} \in R^{+}(w^{-1})$ 
implies $\beta_{i_1,j_1} +\beta_{i_2,j_2}$ is not a root. This follows from 
the fact that no element of $R^{+}(w^{-1})$ starts or ends with 
$\alpha_{a_k+1}$, for any $k$, $1\leq k\leq r-1$ 
(i.e. for all $\beta_{i,j}\in R^{+}(w^{-1})$ and $1\leq k\leq r-1$, 
$\beta_{i,j}-\alpha_{a_k+1} \neq 0$ is not a root.)

Now the natural action of the maximal torus $T$ on $GL_n(k)/P_r$,
induces an action of $T$ on $V$ .

\begin{lemma}
Consider the torus $T^{'}=\prod_{\beta\in R^{+}(w^{-1})} G_{m,\beta}$ 
where $G_{m,\beta}=G_m$ for each $\beta\in R^{+}(w^{-1})$.
We have a natural action of $T$ on $T^{'}$ through the homomorphism
of algebraic groups $\Phi:T\rightarrow T^{'}$ defined by 
$\Phi(t)=( \beta(t))_{\beta}$ for all $t\in  T$. The map
$V\rightarrow T^{'}$ defined by 
$\prod u_{\beta}(x_{\beta})w.P \mapsto (x_{\beta})_{\beta}$
is a $T$-equivariant isomorphism of varieties.
\end{lemma} 
\begin{proof}
Proof is easy.
\end{proof}

We now describe all the Schubert varieties admitting semi-stable points.\\
Let $n=qr+t$, with $1\leq t \leq r$ and let $w\in W^{I_{r}}$.
\begin{lemma}
 Then the following are equivalent:\\
$(1)$ $X(w)_{T}^{ss}(L_{r})$ is non-empty. \\
$(2)$ $\tau_{r}\leq w$ and $w \tau_{r}^{-1}\leq (\tau^{n-r})^{-1}$.\\
$(3)$ $w=(s_{a_{1} }\cdots s_1)\cdots (s_{a_{r}}\cdots s_{r})$, where 
$\{a_i: i=1,2\cdots r\}$ is an increasing sequence of positive integers such 
that $a_i\geq i(q+1)\,\ \forall \,\ i\leq t-1$ and 
$a_i=iq+(t+1)\,\ \forall \,\ t\leq i\leq r$.

\end{lemma} 
\begin{proof}

By Hilbert-Mumford criterion (theorem 2.1 of [3]) a point $x\in G/P_{r}$ is 
semi-stable if and only
if $\mu^{L}(\sigma x, \lambda)\leq 0$ for all $\lambda\in \overline{C(B)}$
and for all $\sigma\in W$. By the lemma 2.1 of [6], this statement is 
equivalent to $\langle -w_{\sigma}(\omega) , \lambda \rangle \geq 0$ for all 
$\lambda \in \overline{C(B)}$ and for all $\sigma\in W$, where 
$w_{\sigma}\in W^{I_{r}}$ is such that $\sigma x \in U_{w_{\sigma}}w_{\sigma}P_{r}$. 
Thus, by corollary 1.8 applied to the situation $\omega = \omega_{r}$, a point 
$x$ is semi-stable if and only if $x$ is not in the W- translates of 
$U_{\tau}\tau P_{r}$ with $\tau\in W^{I_{r}}$ and $\tau_{r}\not \leq \tau$. 

Now, for a $w\in W^{I_{r}}$, $X(w)$ is not contained in the finite union 
$\bigcup_{\tau\not\geq \tau_{r}}U_{\tau}\tau P_{r}$ if and only if $\tau_{r}\leq w$.
The second condition $w\tau_{r}^{-1}\leq (\tau^{n-r})^{-1}$ is an immediate 
consequence when $w \geq \tau_{r}$. 
This completes the proof.
\end{proof}

\begin{proposition}
Let $X_{i,j}$ denote the regular function on $V$ defined by 
$\prod u_{\beta_{kl}}(x_{\beta_{kl}})w.P \mapsto x_{\beta_{ij}}$
for all $1\leq i\leq r-1$ and $1\leq j\leq a_i-i+1$; and
let $Y_{i,j}:=\frac{X_{i,a_i-i+1}.X_{i+1,j}}{X_{i,j}.X_{i+1,a_i-i+1}}$.
Then the ring of $T$-invariant regular functions is generated
by $Y_{i,j},Y_{i,j}^{-1}$, where $1\leq j\leq a_i-i$, for each $i$, 
and $1\leq i\leq r-1$; $Y_{i,j}$ are algebraically independent.
\end{proposition}

\begin{proof}

Now, consider the homomorphism of tori,
\[T\stackrel{\Psi}{\longrightarrow}T'\,\ defined\,\  by\]  
\[ \Psi(t)= (t^{\beta_{ij}}),\,\ i=1,2\cdots r,j=1,2\cdots a_i-i+1.\]
Proof of the proposition follows from the following claim.\\
\underline{Claim}: $E_{i,a_i-(i-1)}-E_{i+1,a_i-(i-1)}-E_{i,j}+E_{i+1,j} ; 
i=1,2\cdots r-1\,\ and\,\ j=1,2\cdots a_i-i$ forms a basis for 
$Ker(\Psi^{*}:X(T){\longrightarrow}X(T'))$, \,\ where\,\ $E_{i,k}$ 
is the matrix with 1 in the $(i,k)^{th}$ place and 0 elsewhere.\\ 
\underline{Proof of the claim}: Now any character of $T^{'}$ is of the form
$(t_{\beta})\mapsto \prod t_{\beta}^{m_{\beta}}$ where $m_{\beta}$ are 
integers. Now such a character is $T$-invariant iff the sum 
$\sum_{\beta} m_{\beta}\beta$ is zero.
Plugging in the expression of $\beta$'s in terms of the simple roots
$\alpha_k$'s and noting that they are linearly independent we get a set
of linear equations over $\mathbb Z$, by equating to zero the coefficient
of each $\alpha_k$. 
Let us denote by $R(p)$, $1\leq p\leq r$ 
the set of roots appearing in $p$-th row of the array described above; 
and let $C(q)$, $1\leq q\leq a_r-(r-1)$ denote the set of roots 
appearing in the $q$-th column of the array.

Comparing the coefficient of $\alpha_1$, we have 
\,\ $\sum_{\beta\in C(1)}m_{\beta}=0$ .\\
Comparing the coefficient of $\alpha_2$, and using the above observation, we get
\,\ $\sum_{\beta\in C(2)}m_{\beta}=0$.
Proceeding this way, we get  
\[\sum_{\beta\in C(j)}m_{\beta}=0\,\,\, \forall j,\,1\leq j\leq a_1.\]
Let $k$ be the least positive integer such that $\alpha_k+\cdots +\alpha_{a_i}$
is the first root in the column $C(a_1+1)$.\\
Comparing the coefficient of $\alpha_k$, we get \,\
$\sum_{\beta\in C(a_1+1)}m_{\beta}=0$.\\
Proceeding this way, we get 
\[ \sum_{\beta\in C(j)}m_{\beta}=0\,\,\, \forall j,\,1\leq j\leq a_r-r+1.\]
Now comparing the coefficient of $\alpha_{a_r}$, we get \,\
$\sum_{\beta\in R(r)}m_{\beta}=0$.\\
Comparing the coefficient of 
$\{\alpha_j: j=a_{r-1},2+a_{r-1},\cdots a_r\}$, we get 
\[\sum_{\beta\in R(r-1)}m_{\beta}+\sum_{\beta\in R(r)}m_{\beta}=0.\]
Thus we have \[\sum_{\beta\in R(r-1)}m_{\beta}=0.\]
Proceeding this way, we get 
\[\sum_{\beta\in R(i)}m_{\beta}=0 \,\,\, \forall i,\,1\leq i\leq r. \]
\end{proof}

\section{Description of the action of the Weyl group on the quotient ${_{_{T}\backslash \backslash}}G_{r,n}^{ss}(\mathcal{L}_r)$ }

In this section, we describe the action of the Weyl group on  the quotient 
${_{_{T}\backslash \backslash}}G_{r,n}^{ss}(\mathcal{L}_r)$.

We first write down the stabiliser of $X(w)$ in $W$. Let $w=(s_{a_1} \cdots s_1)
(s_{a_2}\cdots s_2)\cdots (s_{a_r}\cdots s_r)\in W^{I_{r}}$ be such that 
$w\geq \tau_{r}$. Then, we have

\begin{lemma}
Description of the set $\{s_i: s_i(X(w))\subseteq X(w), i=1,2,\cdots n-1\}$:\\
$1$. $\{s_j:1\leq j\leq a_1-2\}$.\\
$2$. $\{s_j:a_p+2\leq j\leq a_{p+1}-2, p=1,2,\cdots r-1\}$.\\
$3$. $\{s_{a_p-1}: p=1,2,\cdots r\}$.\\
$4$. $\{s_{a_p}: p=1,2,\cdots r\}$. \\
\end{lemma}
\begin{proof}
Proof uses braid relations of the Weyl group $S_n$.
\end{proof}

We now explicitely describe the action of the stabilisers on

\begin{proposition}
\underline{Description of the action}:\\
$1$. $s_j$ interchanges $Y_{i,j}$ and $Y_{i,j+1}$ for $i=1,2,\cdots r-1$, and
keeps all other $Y_{i,k}$'s fixed.\\
$2$. $s_j$ interchanges $Y_{i,j-p}$ and $Y_{i,j-p+1}$ for $p+1\leq i\leq r-1$, 
and keeps all other $Y_{i,k}$'s fixed.\\
$3(a)$. If $2\leq p\leq r$, then $s_{a_p-1}$ fixes all the $Y_{i,k},
1\leq i\leq p-1$.\\
$(b)$. If $p\leq i\leq r-1$, $a_p-p=a_i-i$ and $1\leq k\leq a_p-p$, then
$s_{a_p-1}(Y_{i,a_p-p})=Y_{i,a_p-p}^{-1}$, and $s_{a_p-1}(Y_{i,k})=Y_{i,k}.
Y_{i,a_p-p}^{-1}$.\\
$(c)$. If $p+1\leq i\leq r-1$, $a_i-i\gneq a_p-p$, then $s_{a_p-1}(Y_{i,a_p-p})=
Y_{i,a_p-p+1}$, and keeps all other $Y_{i,k}$'s fixed.\\
$4(a)$.\underline{ $2\leq p\leq r-1$, and $a_p=a_{p-1}+1$}.\\
$(i)$. If $3\leq p\leq r$ and $1\leq k\leq a_{a_{p-2}-p+2}$, then 
$s_{a_p}(Y_{p-2,k})=Y_{p-2,k}.Y_{p-1,k}.Y_{p-1,a_{p-2}-p+3}^{-1}$.\\
$(ii)$. If $1\leq k\leq a_p-p$ then $s_{a_p}(Y_{p-1,k})=Y_{p-1,k}^{-1}$ and 
$s_{a_p}(Y_{p,k})=Y_{p,k}.Y_{p-1,k}$.\\
$(iii)$. $Y_{i,k}$'s are fixed for $i\neq p-2, p-1, p$ and $1\leq k\leq a_i-i$.
\\
$(b)(i)$. If $1\leq i\leq p-1$ or $a_p-p+1\leq k\leq a_r$, $Y_{i,k}$'s are 
fixed.\\
$(ii)$. If $i=p$ and $1\leq k\leq a_p-p$ then $s_{a_p}(Y_{p,k})=1-Y_{p,k}$.\\
$(iii)$ If $p+1\leq i\leq r-1$ and $1\leq k\leq a_p-p$, then, 
$ s_{a_p}(Y_{i,k})= \frac{1-\prod_{m=p}^{i}(Y_{m,k}/Y_{m,a_p-p+1})}
{1-\prod_{m=p}^{i-1}(Y_{m,k}/Y_{m,a_p-p+1})}\times Y_{i,a_p-p+1}$.\\
$(c)$.\underline{ Action of $s_{a_r}$}: \\
$(i)$. If $a_r=a_{r-1}+1$ then $s_{a_r}(Y_{r-2,k})=Y_{r-2,k}.Y_{r-1,k}.Y_{r-1,a_{r-2}-r+3}^{-1}$, for $1\leq k\leq a_{r-2}-r+2$ and $s_{a_r}(Y_{r-1,k})=Y_{r-1,k}^{-1}$, for $1\leq k\leq a_r-r$.\\
$(ii)$. If $a_{r-1}+2\leq a_r$ then $Y_{r,k}$'s are fixed for 
$1\leq k\leq a_r-r+1$.
\end{proposition}

\begin{proof}

Proof is essntially based on the following properties of groups with $BN$-pair 
and commutator relations:

$(i)$
$ \left( \begin{array}{cc}
0 & 1 \\
1 & 0 \\
\end{array} \right)$ 
$ \left( \begin{array}{cc}
1 & x \\
0 & 1 \\
\end{array} \right)$ 
$ \left( \begin{array}{cc}
0 & 1 \\
1 & 0 \\
\end{array} \right)$ =$ \left( \begin{array}{cc}
1 & \frac {1}{x} \\
0 & 1 \\
\end{array} \right)$ $ \left( \begin{array}{cc}
0 & 1 \\
1 & 0 \\
\end{array} \right)$ $ \left( \begin{array}{cc}
x & 1 \\
0 & \frac {-1}{x} \\
\end{array} \right)$ ,  and

$(ii)$ $[u_{\alpha}(x_{\alpha}),u_{\beta}(x_{\beta})] = \left\{ \begin{array}{ll}
         u_{\alpha+\beta}(x_{\alpha}.x_{\beta}) & 
\mbox{if $\alpha=\epsilon_i-\epsilon_j\,\ and\,\ \beta=\epsilon_j-\epsilon_k 
\,\,, i<j<k$};\\
        u_{\alpha+\beta}(-x_{\alpha}.x_{\beta}) & \mbox{if 
$\alpha=\epsilon_i-\epsilon_j\,\ and\,\ \beta=\epsilon_k-\epsilon_i 
\,\,, k<i<j$}.\end{array} \right.$ 

We first consider the action of $W$ on the $X_{j,k}$'s and then describe 
resulting action on the $Y_{j,k}$'s. If $1\leq i\leq a_1-2$
then $s_i$ interchanges $X_{j,i}$ and $X_{j,i+1}$ for all $j$, $1\leq j\leq r$.
Therefore, it follows that $s_i$ interchanges $Y_{j,i}$ and $Y_{j,i+1}$ for all 
$j$, $1\leq j\leq r-1$ and keeps all other $Y_{j,k}$'s fixed.
Similarly for $p\geq 2$ and $a_p+2\leq a_{p+1}$, if 
$a_p+2\leq i\leq a_{p+1}-2$, $s_i$ interchanges $X_{j,i-p}$ and $X_{j,i-p+1}$.
Thus $s_i$ interchanges $Y_{j,i-p}$ and $Y_{j,i-p+1}$ for 
all $j$, $i+1\leq j\leq r-1$ and keeps all other $Y_{j,k}$s fixed. 
Now, we compute the actions of $s_{a_i-1}$, $s_{a_i}$ and $s_{a_i+1}$. 

\underline{Action of $s_{a_i+1}$ for each $i$, $1\leq i\leq r-1$}\\
\underline{{\em Case I}: $a_i+2\leq a_{i+1}$} \hspace{.5cm}In this case we have 
\[ 
\begin{array}{l}
s_{a_i+1}w \\
=s_{a_1+1}.(s_{a_1}\ldots s_1).(s_{a_2}\ldots s_2)\ldots (s_{a_r}\ldots s_r)\\
=(s_{a_1}\ldots s_1)\ldots(s_{a_i+1}.s_{a_i}\ldots s_i)\ldots (s_{a_r}\ldots s_r)
\end{array}
\]
which is a reduced expression and $s_{a_i+1}.w\in W^I$ by lemma $1.12$.
Now lemma $1.13$ implies that $s_{a_i+1}.w\geq w$. Hence, $X(w)$ is not 
stable under the action of $s_{a_i+1}$. \\
\underline{{\em Case II}: $a_i+1=a_{i+1}$} \hspace{.5cm}In this case 
$s_{a_i+1}=s_{a_{i+1}}$ 
and the
action will be described in the later part of this paragraph. In fact
we see that in this case $(s_{a_i+1}w)^I=w$. Hence $X(w)$ is stable under
the action of $s_{a_i+1}$. 

\underline{Action of $s_{a_i-1}$}\\
In case $i=1$, we may assume that $a_1\neq 1$, and for $i\geq 2$, 
$a_{i-1}\neq a_i-1$. Now $s_{a_i-1}$ interchanges the
$(a_i-i)$-th and $(a_i-i+1)$-th columns of each of the $j$-th row, 
of the array of roots $R^{+}(w^{-1})$, for
$i\leq j\leq r$; thus $s_{a_i-1}$ interchanges $X_{j,a_i-i}$ and $X_{j,a_i-i+1}$
for each $j$, $i\leq j\leq r$. Therefore, the action of $s_{a_i-1}$ is as 
follows:\\
$(1)$ $s_{a_i-1}$ fixes all the $Y_{j,k}$, for $1\leq j\leq i-1$, for 
$i\geq 2$. \\
$(2)$ For $j\geq i\leq r-1$ and $a_i-i=a_j-j$, 
$Y_{j,a_i-i}\mapsto Y_{j,a_i-i}^{-1}$, and 
for $Y_{j,k}\mapsto Y_{j,k}.Y_{j,a_i-i}^{-1}$ for $1\leq k<a_i-i$.\\
$(3)$ For $i+1\leq j\leq r-1$ if $a_j-j>a_i-i$, then $s_{a_i-1}$
interchanges $Y_{j,a_i-i}$ and $Y_{j,a_i-i+1}$ and keeps all other
$Y_{j,k}$'s fixed.

\underline{ Action of $s_{a_i}$ for $1\leq i\leq r$}\\
Let us show that $X(w)$ is stable under the action of each of
the $s_{a_i}$.
Let 
\[w=(s_{a_1}\ldots s_1).(s_{a_2}\ldots s_2)\ldots (s_{a_r}\ldots s_r)\]
Thus 
\[s_{a_i}w=(s_{a_1}\ldots s_1)\ldots (s_{a_{i-2}}\ldots s_{i-2}).s_{a_i}.
(s_{a_{i-1}}\ldots s_{i-1}).(s_{a_i}\ldots s_i)\ldots (s_{a_r}\ldots s_r)\]
\emph{ Case $1$}: $i=1$, or $a_{i-1}+2\leq a_i$ for $i\geq 2$. 
In this case it is clear that
\[s_{a_i}w=(s_{a_1}\ldots s_1)\ldots (s_{a_{i-2}}\ldots s_{i-2}).
(s_{a_{i-1}}\ldots s_{i-1}).(s_{a_i-1}\ldots s_i)\ldots (s_{a_r}\ldots s_r)\]
which, by lemma $1.12$ and $1.13$, is in $W^I$ and $s_{a_i}w\leq w$. \\
\emph{Case $2$}: $a_{i-1}+1=a_i$. Note that,
\[ w_1=(s_{a_{i-1}}\ldots s_{i-1}).(s_{a_i}\ldots s_i)\in W^J\]
where $J=S\setminus \{\alpha_i\}$. Now, 
\[
\begin{array}{rcl}
w_1(\omega_i)& = &\omega_i-\sum_{j=i-1}^{a_{i-1}}\alpha_j-\sum_{j=i}^{a_i}\alpha_j 
\\
\Rightarrow s_{a_i}w_1(\omega_i)& = &s_{a_i}(\omega_i)-\sum_{j=i-1}^{a_{i-1}}
\alpha_j-\sum_{j=i}^{a_i}\alpha_j 
\end{array}
\]
Now, if $a_i=i$, then $a_{i-1}=i-1$; so 
$s_{a_i}w_1=s_i.s_{i-1}.s_i=s_{i-1}.s_i.s_{i-1}=w_1.s_{i-1}$. Otherwise,
$a_i\neq i$. This implies that $s_{a_i}(\omega_i)=\omega_i$. Therefore,
$s_{a_i}w_1(\omega_i)=w_1(\omega_i)$. Hence, by lemma $1.3$, we get
$s_{a_i}w_1=w_1.s_{\alpha}$ for some $\alpha\in J$. This gives 
$w_1^{-1}s_{a_i}w_1=s_{w_1^{-1}(\alpha_{a_i})}=s_{\alpha}$. Now it follows 
that $w_1^{-1}(\alpha_{a_i})=\alpha_{i-1}$. Hence, $s_{a_i}w_1=w_1.s_{i-1}$.
Therefore, in both the sub-cases $s_{a_i}.w=w.s_{i-1}$; 
in particular $(s_{a_i}.w)^I=w$.
Now we shall compute the action of $s_{a_i}$,
for $1\leq i\leq r$.

\underline{\emph{ Case I:} $2\leq i\leq r-1$ and $a_{i}=a_{i-1}+1$}.
In this case,
$s_{a_i}$ interchanges $X_{i,k}$ and $X_{i-1,k}$ for $1\leq k\leq a_i-i+1$
and keeps all other $X_{j,k}$'s fixed. Hence, the action of $s_{a_i}$ on
the $Y_{j,k}$'s is as follows:\\ $(1)$ If $i\geq 3$, 
$Y_{i-2,k}\mapsto Y_{i-2,k}.Y_{i-1,k}.Y_{i-1,a_{i-2}-i+3}^{-1}$ for 
$1\leq k\leq a_{i-2}-i+2$ \\ $(2)$ $Y_{i-1,k}\mapsto Y_{i-1.k}^{-1}$ for
$1\leq k\leq a_i-i$. \\ $(3)$ $Y_{i,k}\mapsto Y_{i,k}.Y_{i-1,k}$ for 
$1\leq k\leq a_i-i$. \\ $(4)$ $Y_{j,k}$ is fixed for $1\leq k\leq a_j-j$
for each $j\neq i-2,i-1,i$.

\underline{\emph{Case II:} $a_i\geq a_{i-1}+2$ for $2\leq i\leq r-1$, 
or $i=1.$} 
In this case 
$s_{a_i}$ changes only the $i$-th row and the $(a_i-i+1)$-th column of the
array of roots $R^{+}(w^{-1})$. The resulting $i$-th row turns out to be
\[\alpha_1+\alpha_2+\cdots +\alpha_{a_i-1},\,\alpha_2+\cdots +\alpha_{a_i-1},\,
\cdots ,\alpha_{a_1}+\cdots +\alpha_{a_i-1},\,\alpha_{a_1+2}+\cdots +
\alpha_{a_i-1},\cdots \, ,\alpha_{a_i-1},\,\, -\alpha_{a_i}\]
and the transpose of the $(a_i-i+1)$-th column turns out to be
\[-\alpha_{a_i},\, \alpha_{a_i+1}+\cdots+\alpha_{a_{i+1}},\,\alpha_{a_i+1}+\cdots+
\alpha_{a_{i+2}},\,\cdots, \alpha_{a_i+1}+\cdots+\alpha_{a_r}\] 
Let $\beta_{j,k}$ be any root which is fixed under the action of $s_{a_i}$.
and let $\beta_{p,q}$ be any root of the $i$-th row or the $(a_i-i+1)$-th
column, i.e. either $p=i$ or $q=a_i-i+1$. We claim that $u_{\beta_{i,j}}(X_{i,j})$
and $u_{s_{a_i}\beta_{p,q}}(X_{p,q})$ commute. This follows from the fact that
$\beta_{j,k}-\alpha_{a_i}\not \in R^{+}(w^{-1})$ and the observation that
for any root $\beta\in R^{+}(w^{-1})$ and $1\leq m\leq r$ 
$\beta -\alpha_{a_m+1}\not \in R^{+}$.
Let us denote by $M$ the sub-array consisting of $\beta_{k,l}$ where 
$k\geq i$ and $1\leq l\leq a_i-i+1$.
Then
\[
\begin{array}{l}
s_{a_i}.(u_{\beta_{r,1}}(X_{r,1}).u_{\beta_{r,2}}(X_{r,2})\ldots 
u_{\beta_{r,a_r-r+1}}(X_{r,a_r-r+1}).u_{\beta_{r-1,1}}(X_{r-1,1}).u_{\beta_{r-1,2}}
(X_{r-1,2}) \\
\ldots u_{\beta_{r-1,a_{r-1}-r+2}}(X_{r-1,a_{r-1}-r+2})\ldots 
u_{\beta_{1,1}}(X_{1,1}).u_{\beta_{1,2}}(X_{1,2})\ldots 
u_{\beta_{1,a_1}}(X_{1,a_1})).w.P \\
=(\prod_{\beta_{k,l}\not \in M}u_{\beta_{k,l}}(X_{k,l})).s_{a_i}.(u_{\beta_{r,1}}
(X_{r,1}).u_{\beta_{r,2}}(X_{r,2})\ldots u_{\beta_{r,a_i-i+1}}(X_{r,a_i-i+1})).
(u_{\beta_{r-1,1}}
(X_{r-1,1})\cdots \\
u_{\beta_{r-1,2}}(X_{r-1,2})\ldots u_{\beta_{r-1,a_i-i+1}}(X_{r-1,a_i-i+1}))\ldots
u_{\beta_{i,1}}(X_{i,1}).u_{\beta_{i,2}}(X_{i,2})\ldots u_{\beta_{i,a_1-i+1}}
(X_{i,a_i-i+1})).w.P
\end{array}
\]

Thus the action of $s_{a_i}$, in this case is as follows:
\[X_{i,a_i-i+1}\mapsto  X_{i,a_i-i+1}^{-1}; \,\,\,
X_{i,k}\mapsto X_{i,k}.X_{i,a_i-i+1}^{-1} \,\,\mbox{for} \,\,k\leq a_i-i\] 
\[X_{j,k}\mapsto X_{j,k}-\frac{X_{j,a_i-i+1}.X_{i,k}}{X_{i,a_i-i+1}}\,\,
\mbox{ for}\,\, i+1\leq j\leq r\,\,\mbox{and}\,\, 1\leq k\leq a_i-i\]
\[ X_{j,a_i-i+1}\mapsto -X_{j,a_i-i+1}/X_{i,a_i-i+1}\,\,\mbox{for}\,\,
i+1\leq j\leq r\] 

From this the resulting action on the $Y_{j,k}$ turns out to be as 
follows:\\
$(1)$ $s_{a_i}$ fixes $Y_{j,k}$'s provided $j\leq i-1$ or $k\geq a_i-i+1$.\\  
We now make the convention that $Y_{j,k}:=1$ if $k\geq a_j-j+1$ or if 
$j\geq r$. \\
$(2)$ $j=i$. Here, for $k\leq a_i-i$,
\[ 
\begin{array}{rrcl}
&Y_{i,k}&=&\frac{X_{i,a_i-i+1}.X_{i+1,k}}{X_{i+1,a_i-i+1}.X_{i,k}}\\
\therefore & s_{a_i}(Y_{i,k})&=&\frac{X_{i,a_i-i+1}^{-1}.
(X_{i+1,k}-\frac{X_{i+1,a_i-i+1}.X_{i,k}}{X_{i,a_i-i+1}})}
{X_{i,k}.X_{i,a_i-i+1}^{-1}.(-X_{i+1,a_i-i+1}
/X_{i,a_i-i+1})}\\
&&=& 1-Y_{i,k}
\end{array}
\]

$(3)$ $i+1\leq j\leq r-1$ and $1\leq k\leq a_i-i$. 
Define $Y_{j,k}^{'}=(X_{i,a_i-i+1}.X_{j,k})/(X_{j,a_i-i+1}.X_{i,k})$.
Then, we have $s_{a_i}(Y_{j,k})=1-Y_{j,k}$. 
It follows that $Y_{j,k}=Y_{j+1,k}^{'}.{Y_{j,k}^{'}}^{-1}.Y_{j,a_i-i+1}$.
Hence, $s_{a_i}(Y_{j,k})=\frac{1- Y_{j+1,k}^{'}}{1-Y_{j,k}^{'}}.Y_{j,a_i-i+1}$.
Now,
\[
\begin{array}{rcl}
Y_{j,k}^{'}&=&\prod_{m=i}^{j-1}\frac{X_{m,a_i-i+1}.X_{m+1,k}}
{X_{m+1,a_i-i+1}.X_{m,k}}\\
&=&\prod_{m=i}^{j-1}\{(\frac{X_{m,a_m-m+1}.X_{m+1,k}}{X_{m+1,a_m-m+1}.X_{m,k}})\times
(\frac{X_{m,a_m-m+1}.X_{m+1,a_i-i+1}}{X_{m+1,a_m-m+1}.X_{m,a_i-i+1}})^{-1}\}\\
&=&\prod_{m=i}^{j-1}(Y_{m,k}/Y_{m,a_i-i+1})
\end{array}
\]
Thus we have,
\[
\mbox{
\framebox{$
 s_{a_i}(Y_{j,k})= \frac{1-\prod_{m=i}^{j}(Y_{m,k}/Y_{m,a_i-i+1})}
{1-\prod_{m=i}^{j-1}(Y_{m,k}/Y_{m,a_i-i+1})}\times Y_{j,a_i-i+1}$ } }
\]
\emph{ Case III}: \underline{ Action of $s_{a_r}$}: 
$(1)$ If $a_r=a_{r-1}+1$, then $s_{a_r}$ interchanges 
$X_{r-1,k}$ and $X_{r,k}$, $1\leq k\leq a_r-r+1$. A straightforward
checking proves as in \emph{Case I} above, that in this case the action 
of $s_{a_r}$ is as follows:
\[
\begin{array}{ll}
Y_{r-2,k}\mapsto Y_{r-2,k}.Y_{r-1,k}.Y_{r-1,a_{r-2}-r+3}^{-1} &\,\mbox{for}\,
1\leq k\leq a_{r-2}-r+2\\  
Y_{r-1,k}\mapsto Y_{r-1.k}^{-1} &\,\mbox{for}\,
1\leq k\leq a_r-r
\end{array}
\]
$(2)$ If $a_r\geq a_{r-1}+2$, $s_{a_r}$ changes only
$X_{r,k}$'s for $1\leq k\leq a_r-r+1$, as follows:\\
\[
\begin{array}{ll}
X_{r,k}\mapsto X_{r,k}.X_{r,a_r-r+1}^{-1}& \mbox{for}\, 1\leq k\leq a_r-r\\
X_{r,a_r-r+1}\mapsto X_{r,a_r-r+1}^{-1} & \\
\end{array}
\]
 It can be easily checked from here that the $Y_{i,j}$'s are all fixed by 
$s_{a_{r}}$.
\end{proof}

\section{A stratification of ${_{_{N}\backslash \backslash}}
G_{2,n}^{ss}(\mathcal{L}_2)$.}

In this section, we give a stratification of ${_{_{N}\backslash \backslash}}
G_{2,n}^{ss}(\mathcal{L}_2)$.
\begin{lemma}
Let $w\in W^{I_2}$. Let $x\in{U_wwP_2}^{ss}$ be such that $x$ is not in the $W$
-translate of $X_{\tau}, \tau<w$. If $\sigma(x)\in {U_wwP_2}$, then $\sigma 
\in$ Stabiliser of $X(w)$ in $W$.
\end{lemma}

\begin{proof}
Let $\sigma \in W$ be of minimal length such that $\sigma x \in {U_wwP_2}$.
Then $\sigma=\sigma_1.\sigma_2$ with 
$l(\sigma)=l(\sigma_1)+l(\sigma_2)$ and $\sigma_2.w \in W^I, w\leq \sigma_2w$.\\
Let $\sigma_2$ be of maximal length with this property.
So $\sigma.w=s_{m+t+1}s_{m+t}\cdots s_{m+1}w, t\geq 1$, and 
$w=(s_m\cdots s_1)(s_{n-1}\cdots s_2)$.\\
Now, $\sigma_1(\sigma_2u_wwP_2)\in U_wwP_2$. $\hspace{9cm}\ldots (1)$\\
Since $\sigma_2$ is of maximal length $s_{m+j}\nleq \sigma_1$ for some 
$j\geq 1$. $\hspace{4cm}\ldots (2)$\\
Now, $\sigma_2x\in U_{{\sigma}_{2}w}\sigma_2wP_2$. Since 
$l(\sigma)=l(\sigma_1)+l(\sigma_2)$ and $\sigma^{-1}(\alpha_{m+t+1})<0$, 
$\sigma_2$ is of maximal length, we may assume that $\sigma_1(\alpha_j)>0$.  
$\hspace{7cm}\ldots (3)$  \\
From $(1),(2)$ and $(3)$, $\sigma_1$ must take a reduced form as 
\[
\begin{array}{rcl}
\sigma_1&=&(\phi s_{m+t-1}s_{m+t+1}s_{m+t})\sigma_2\\
&=&\phi s_{m+t-1}(s_{m+t+1}s_{m+t}s_{m+t+1})s_{m+t}\sigma_{2}^{\prime}\\
&=&\phi s_{m+t-1}s_{m+t}s_{m+t+1}\sigma_2^{\prime}
\end{array}
\]
This contradicts the assumption that $l(\sigma)=l(\sigma_1)+l(\sigma_2)$.\\
This completes the proof.
\end{proof}

The longest element of $W^{I_2}$ is 
\[w^I_0=(s_{n-2}.s_{n-3}\ldots s_1).(s_{n-1}.s_{n-2}\ldots s_2)\] and the 
unique minimal element $\tau_2$ of $W^I$ such that
$\tau_{2}(n\omega_2)\leq 0$ is
\[\tau_{2}=(s_{\lceil \frac{n-1}{2}\rceil}.s_{\lceil \frac{n-1}{2}\rceil-1}\ldots 
s_1).(s_{n-1}.s_{n-2}\ldots s_2)\]
Therefore any element $w\in W^I$ such that $X(w)_T^{ss}(\mathcal{L}_2)\neq 
\empty$ is of the 
form 
\[ w=(s_m.s_{m-1}\ldots s_{\lceil \frac{n-1}{2}\rceil}.s_{\lceil \frac{n-1}
{2}\rceil-1}\ldots s_1).(s_{n-1}.s_{n-2}\ldots s_2)\] with $m\geq \lceil
\frac{n-1}{2}\rceil$.

\begin{proposition}
Let $r=2, w=(s_m...s_1)(s_{n-1}...s_2),\lceil{\frac{n-1}{2}}\rceil\leq m\leq 
n-2$. We can arrange the $Y_{ij}$'s as $Y_1,Y_2,\cdots Y_{m-1}$ with 
\[ s_i(Y_i)=Y_{i+1},\] \[s_i(Y_j)=Y_j \,\ if \,\ j=k,i+1\,\ and\,\ i=1,2\cdots
 m-2,\]
\[s_{m-1}(Y_i)=Y_i.Y_{m-1}^{-1},\,\ if \,\ i\leq m-2,\]\[s_{m-1}(Y_{m-1})=
Y_{m-1}^{-1},\] \[s_m(Y_i)=1-Y_i \,\ for \,\ i=1,2,\cdots m-1.\]
Further, we have 
\[s_i(Y_j)=Y_j ~  \forall ~ i=m+2,\cdots n-1, ~ when ~ m\leq n-3\]
and \[s_{n-1}(Y_j)=Y_{j}^{-1} ~ \forall ~ j ~ when ~ m=n-2.\]
\end{proposition}
\begin{proof}
Proof follows from the proposition 3.2.
\end{proof}

Let $w$ be as in the proposition 4.2.
Now, let $T_{m-1}$ be a maximal torus of $\mathbb PGL_m$, $R_m$ is the root 
system of   $\mathbb PGL_m$. Here, the Weyl group is $S_{m}$, the symmetric group on $m$ symbols.\\
Let $U=\{t\in T:e^{\alpha}(t)\neq 1, \,\ \alpha\in R_m\}$.
Clearly, $U$ is $S_{m}$-stable. On the other hand, $S_{m}$ stabilises 
$(U_{w}wP_{2}/P_{2})_T^{ss}(\mathcal{L}_2)$. 
Let $Y(w)={_{_{T}\backslash \backslash}}(U_{w}wP_{2})_T^{ss}(\mathcal{L}_2)$. Then,
we have 

\begin{corollary}
There is a $S_{m}$-equivariant isomorphism 
$\Psi_{1}:Y(w) \stackrel \sim \longrightarrow U$ such that \\
$\Psi_1^*(e^{\alpha_i+\cdots +\alpha_{m-1}})=Y_i$, \,\ $1\leq i\leq m-1$.
\end{corollary}
\begin{proof}
Proof follows from proposition 4.2.
\end{proof}

Let $\mathfrak{h}_{m}$ be a Cartan subalgebra of 
$\mathfrak{sl}_{m+1}$, $\mathbb P(\mathfrak{h}_m)$ be the projective space and 
$R_{m}\subseteq \mathfrak{h}_{m}^{*}$ be the root system. Let $V_m$ be the open 
subset of  $\mathbb P(\mathfrak{h}_{m})$ defined by 
\[V_m:=\{x\in \mathbb P(\mathfrak{h}_{m}):\alpha(x)\neq 0, \forall \alpha 
\in R_{m}\}.\]
Clearly $V_m$ is $S_{m+1}$-stable.
\begin{corollary}
Let $w=(s_m...s_1)(s_{n-1}...s_2),\lceil{\frac{n-1}{2}}\rceil\leq m\leq 
n-2$.Then, there is a $S_{m+1}$-equivariant isomorphism $\Psi_2:Y(w) 
\stackrel \sim \longrightarrow V$ of affine varieties.
\end{corollary}
\begin{proof}
For $i=1,2\cdots m-1$, take $Z_i=\frac {\alpha_i+\cdots +\alpha_m}{\alpha_m}$ 
and define $\Psi_2$ such that $\Psi_2^*(Z_i)=Y_i$.  
\end{proof}

With notations as above and taking $t=[\frac{n-1}{2}]$ and 
$m= \lceil{\frac{n-1}{2}}\rceil$ we have \\
{\bf Theorem:}
${_{_{N}\backslash \backslash}}G_{2,n}^{ss}(\mathcal{L}_2)$  has a stratification 
$\bigcup_{i=0}^{t}C_i$ where $C_0= {_{S_{m+1}}\backslash}{\mathbb P(\mathfrak{h}_m)}
$, and $C_i= {_{S_{i+m+1}\backslash}{V_{i+m}}}$.
\begin{proof}
Proof follows from lemma 4.1, proposition 4.2 and corollary 4.4.
\end{proof}

\section{Flag variety as a GIT quotient of flag variety of higher dimension}
Let $G=GL_{n+1}(k)$. Let $T$ be a maximal torus of $SL_{n+1}(k)$. Let $B_{n+1}$ 
be a Borel subgroup of $G$ containing $T$. Let 
$S=\{\alpha_{i}:i=1,2,\cdots n\}$ denote 
the set of simple roots with respect to $B_{n+1}$, let $W=S_{n+1}$ be the Weyl 
group. Let $s_{i}$ be the simple reflection corresponding to the simple root 
$\alpha_{i}$. Let $I:=S\setminus \{\alpha_{n}\}$, let $W_{I}$ be subgroup of 
$W$ generated by $\{s_{i}: i \in I\}$ and $w_{0,I}$ denote the longest element 
of $W_{I}$.

\begin{lemma}
Let $\chi=\sum_{i=1}^{n}m_i\alpha_i$ be a regular dominant character, where 
$m_i\in\mathbb N$, $m_{i+1}>m_i$ for $1\leq i\leq n-1$.
Let $w\in W$. Then $w(\chi)\leq 0 \Leftrightarrow w=s_1.s_2\ldots s_n.\tau$ 
for some $\tau\in W_I$.
\end{lemma}
\begin{proof}
$\Rightarrow:$ Since $\chi$ is dominant and $\tau\leq w_{0,I}$, for all
$\tau\in W_I$, we have $\tau(\chi)\geq w_{0,I}(\chi)$; using the fact
that $w_{0,I}(\alpha_i)=-\alpha_{n-i}$ for $i=1,\cdots n-1$ and 
$w_{0,I}(\alpha_n)=\alpha_1+\alpha_2+\cdots +\alpha_n$ we have
$w_{0,I}(\chi)=\sum_{i=1}^{n-1}(m_n-m_{n-i})\alpha_i+m_n.\alpha_n$.
Therefore, $\tau(\chi)=\sum_{i=1}^{n-1}a_i\alpha_i+m_n\alpha_n$, $a_i>0$.
Now, let $w=\phi\tau$ with $\phi\in W^I$, $\tau\in W_I$. Therefore,
$w(\chi)=\phi(\tau(\chi))=\phi(\sum_{i=1}^{n-1}a_i\alpha_i+m_n\alpha_n)$.
Thus $w(\chi)\leq 0$ implies that $\phi=s_1.s_2\ldots s_n$.\\
$\Leftarrow:$ Let $w=s_1.s_2\ldots s_n.\tau$, $\tau\in W_I$. Now,
\[
\begin{array}{rcl}
w(\chi) & = & s_1.s_2\ldots s_n\tau(\chi)\\
        & = & s_1.s_2\ldots s_n(\sum_{i=1}^{n-1}a_i\alpha_i+m_n\alpha_n)\\
        & = & -m_n\alpha_1+\sum_{i=2}^{n}(a_{i-1}-m_n)\alpha_i
\end{array} \]
Since $\chi$ is a dominant weight we have $\chi-\tau(\chi)\geq 0$. Hence
we have $a_i\leq m_i\leq m_n$. Thus $w(\chi)\leq 0$. This completes the proof.
\end{proof}

Consider $GL_{n}(k)$ as a subgroup of $GL_{n+1}(k)$ given by the inclusion $g \mapsto  \left( \begin{array}{cc}
g & 0 \\
0 & 1 \\
\end{array} \right)$. Let $B_{n}=B_{n+1}\bigcap GL_{n}(k)$ as a Borel subgroup with $I$ as the simple roots.

Let $\chi$ be a regular dominant character as in Lemma (5.1).

\begin{theorem}
We have an isomorphism \[ \Psi: {_{_{T}\backslash \backslash}}({GL_{n+1}(k)/B_{n+1}}))^{ss}(L_{\chi})\stackrel \sim \longrightarrow  
GL_n(k)/B_n.\]
\end{theorem}

\begin{proof}

Proof uses cellular decomposition of both homogeneous spaces $GL_{n+1}(k)/B_{n+1}$ 
and $GL_{n}(k)/B_{n}$. First, we fix a total order on the set of positive roots of
$B_{n+1}$ such that $\sum_{i=1}^{n}\alpha_{i} > \sum_{i=1}^{n-1}\alpha_{i} > 
\cdots > \alpha_1 > \sum_{i=2}^{n}\alpha_{i} > \cdots > \alpha_2 > \sum_{i=3}^{n}
\alpha_i > \cdots > \alpha_3 > \cdots >\alpha_{n-1}+\alpha_n > \alpha_n$.  
Now any $GL_{n+1}/B_{n+1}$ (resp. $GL_{n}/B_{n}$) is the union of cells 
$U_{w}wB_{n+1}$ (resp. $U_{\tau}\tau B_{n}$) with $w\in W$ (resp. 
$\tau\in W_{I}$). Using the total order above we can write each element 
$x\in U_{w}$ as a product of $u_{\alpha}$ in the decreasing order from 
the left to the right. Let $X_{\alpha}$ (resp. $Y_{\beta}$) be the co-ordinate 
function on $U_{w}wB_{n+1}$ (resp. $U_{\tau}\tau B_{n}$) corresponding to the 
root $\alpha$ (resp. $\beta$).

With these notations we proceed the proof:

Let $\tau\in W_{I}$. Let $w:= s_{1} s_{2}\cdots s_{n} \tau$.
$V^{0}_{\tau}:\{x=U_{w} wB_{n+1}:X_{\alpha}(x)\neq 0 ~ \forall ~ \alpha\geq 
\alpha_{1}\}$.

Set $V^0:=\bigcup_{\tau\in W_{I}}V_{\tau}^{0}$.

Step1: We prove that $(GL_{n+1}(k)/B_{n+1})^{ss}(L_{\chi})\subset V^0$.

This can be seen as follows.

By Hilbert-Mumford criterion [see theorem 2.1 of [3]], a point 
$x\in GL_{n+1}/B_{n+1}$ 
is semi-stable $\Leftrightarrow$ $\mu^L(x,\lambda)\geq 0$ for all $1$-
parameter subgroup $\lambda$ of $T$ 
$\Leftrightarrow$ $\mu^L(\sigma x,\lambda)\geq 0$ for all one parameter 
subgroups $\lambda \in \overline{C(B)}$ and for all $\sigma \in W$.
By the lemma 2.1 of [6], this statement is equivalent to $\langle -w_{\sigma}
\chi,\lambda\rangle \geq 0$ for all $\lambda \in \overline{C(B)}$ where 
$\sigma x\in U_{w_{\sigma}}w_{\sigma}B$. But this is equivalent to 
$w_{\sigma}(\chi)\leq 0$. And this is equivalent to $w_{\sigma}$ 
is of the form $(s_1\ldots s_n).\tau_1$
for some $\tau_1\in W_I$. \\
Now let 
$x\in U_{w} wB_{n+1}$ with $w=(s_1\ldots s_n)\tau$, $\tau\in W_I$. 

Now, let $X_{\alpha}(x)=0$ for some $\alpha\geq \alpha_{1}$. Let 
$\alpha=\sum_{j=1,\cdots i}\alpha_{j}$. Then, we have $s_{1}s_{2}\cdots s_{i}
x=u^{\prime}\phi B_{n+1}$ with $\phi\neq s_{1}\cdots s_{n}\tau $ for any 
$\tau\in W_{I}$. Hence, by the above discussion, $x$ is not semi-stable.

Step 2: $(GL_{n+1}(k)/B_{n+1})^{ss}(L_{\chi})=V^{0}$. This can be seen by the
above discussion and from the following {\it claim}.

{\it claim:} $V$ is $W$-stable.

\underline{Proof of {\it claim}}: Let $\tau\in W_{I}$.
Let $x \in U_{s_{1}s_{2}\cdots s_{n}\tau}s_{1}s_{2}\cdots s_{n}\tau B_{n+1}$,
with $X_{\alpha}(x)\neq 0$ for all $\alpha\geq \alpha_{1}$.
Then, we have $s_{1}x \in  U_{s_{1}s_{2}\cdots s_{n}\tau}B_{n+1}$ with 
$X_{\alpha}(s_1x)=-\frac{X_{\alpha}(x)}{X_{\alpha_1}(x)}$ for 
$\alpha > \alpha_{1}$, and  $X_{\alpha_{1}}(s_1x)=\frac{1}{X_{\alpha_{1}}(x)}$. 
Hence, $s_{1} x\in V^{0}$.

Now, let $i\neq 1$. If $X_{\alpha_{i}}(u)=0$, then, $s_{i} x= u^{\prime} s_{1} 
s_{2}\cdots s_{n} s_{i-1}\tau B_{n+1}$ with $X_{\alpha}(s_1x)=
X_{s_{i}(\alpha)}(x)$. Hence, $s_{i}(x)\in V^{0}$. Otherwise, we must have
$s_{i} x \in U_{s_{1}s_{2}\cdots s_{n}\tau}B_{n+1} $ with 
$X_{\alpha}(s_ix)=X_{\alpha}(x)$ for all such that $s_{i}(\alpha)=\alpha$,
$X_{\alpha}(s_ix)= \frac{X_{\alpha}(x)}{X_{\alpha_{i}}(x)}$ for all $\alpha$
of the form $\alpha=\sum_{j=k}^{i}\alpha_{j}$ such that $k < i$, 
$X_{\alpha_{i}}(s_ix) =\frac{1}{X_{\alpha_{i}}(x)}$, and $X_{\alpha}(s_ix)=
\frac{-X_{\alpha}(x)}{X_{\alpha_{i}}(x)}$ for all $\alpha$ of the form 
$\alpha= \sum_{j=i}^{k}\alpha_{j}$  such that $k >i$.   

Hence $s_{i} V^{0}\subset V^{0}$ for all $i=1, \cdots n$. Thus, the {\it claim} 
follows from the fact that $W$ is generated by $s_{i}$'s. 
  
Step 3: Now, for each $\tau\in W_{I}$, we exhibit an isomorphism 
\[\Psi_{\tau}:_{_{T\backslash \backslash}}V_{\tau}^{0}\stackrel \sim \longrightarrow 
U_{\tau}\tau B_{n}/B_{n}.\]

Let $\tau\in W_{I}$, consider the map  
$\pi_\tau: V_{\tau}^{0}\longrightarrow (U_{\tau}\tau B_n)/B_n$
defined by 
$\phi_\tau(x)= y $ with for each $\beta\not\geq \alpha_{1}$ $Y_{s_n\ldots 
s_1(\beta)}(y)=(\frac{-X_{\beta}(x)X_{\beta^{'}}(x)}{X_{\beta+\beta^{'}}(x)})$
where for each $\beta\in R^{+}(w^{-1})$ with $\beta\not\geq \alpha_1$,
$\beta^{'}$ is the unique element of $R^{+}$ with $\beta^{'}\geq \alpha_1$
such that $\beta+\beta^{'}\in R^{+}$. 
Clearly this map is $T$-invariant.
Thus the morphism $\pi_\tau$ give rise to a morphism 
\[\Psi_{\tau}:_{_{T\backslash \backslash}}V_{\tau}^{0}\longrightarrow 
U_{\tau}\tau B_{n}/B_{n}.\]

 Clearly $\Psi_{\tau}$is surjective. We now prove that $\Psi_{\tau}$ is 
injective:

$\pi_w$ is injective for each
$w\in W$ of the form $w=s_1.s_2\ldots s_n\tau$, for some $\tau\in W_I$.
Let $x_1$ and $x_2$ be two points of $V_{\tau}^{0}$ such that 

$\pi_{\tau}(x_1)=\pi_{\tau}(x_2)$. 
Hence, $\frac{X_{\beta}(x_1)X_{\beta^{'}}(x_1)}{X_{\beta+\beta^{'}}(x_1)}=
\frac{X_{\beta}(x_2)X_{\beta^{'}}(x_2)}{X_{\beta+\beta^{'}}(x_2)}$. Let $t\in T$ 
be such that 
$(\alpha_1+\cdots +\alpha_i)(t)=\frac{X_{\alpha_1+\cdots +\alpha_i}(x_2)}
{X_{\alpha_1+\cdots +\alpha_i}(x_1)}$ for all $i$, $1\leq i\leq n$. Then, it is 
easy to check that $t\cdot x=y$. \\
Thus $\Psi_{\tau}$ is bijective for each $\tau\in W_{I}$. 

Step 4: $\Psi_{\tau}$ puts together to give an isomorphism 

\[\Psi: _{_{T\backslash \backslash}}V^{0}\stackrel{\sim}{\longrightarrow} GL_{n}(k)/B_{n}.\]

Since the $W$- translates of $V_{w_{0,I}}^{0}$ is the whole of $V^0$, and 
$W_{I}$- translates of $U_{w_{0,I}}w_{0,I}B_n$ is the whole of $GL_{n}/B_{n}$, 
and there is an isomorphism  from $W_{S\setminus \{{\alpha_1}\}}$ to $W_{I}$ taking
$s_{i}$ to $s_{i-1}$ for each $i=2, \cdots n$, to prove the Theorem, it is 
sufficient to prove that the $T$ - invariant morphisms 
$\pi_{\tau}: V_{\tau}^{0}\longrightarrow U_{\tau}\tau B_n$, and 
$\pi_{\tilde{s_{i-1}\tau}} U_{\tilde{s_{i-1}\tau}}\longrightarrow 
U_{\tilde{s_{i-1}\tau}}$ satisfy the following:

$Y_{\alpha}(\pi_{\tau}(x))=Y_{\alpha}(s_{i-1}(\pi_{\tilde{s_{i-1}\tau}}(s_i x)))$
for each $\alpha\in R^{+}(\tau^{-1})$. (Here, the notation 
$\tilde{s_{i-1}\tau}=\tau$ if $s_{i-1}\tau < \tau$, and $\tilde{s_{i-1}\tau}=
s_{i-1}\tau$ otherwise.)

We make use of the following observations using commutator relations: 

\beqn
X_{\alpha}(s_i; x)= \left\{
\begin{array}{ll}
 \frac{-X_{\alpha}(x)}{X_{\alpha_i}(x)} & \mbox{ if } 
 \alpha = \alpha_i + \cdots \alpha_k,  ~ i < k ~  \mbox{ and } ~ 
w^{-1}(\alpha_{i+1}+\cdots \alpha_k) >0, \\
\frac{1}{X_{\alpha_i}(x)} & \mbox{ if } \alpha=\alpha_i, \\
X_{s_{i}(\alpha)}(x) & otherwise  \end{array}
\right. .
\eeqn

Let $\alpha\in R^{+}(\tau^{-1})$

Case 1:  $\alpha=\alpha_{k-1} + \cdots \alpha_{i-1}$, $k < i$, 
$w^{-1}(\alpha_{k}+\cdots \alpha_{i})=\tau^{-1}(\alpha_{k-1}+\cdots 
\alpha_{i-1}) > 0$ and $s_{i-1}\tilde{\tau}=\tau$. 

In this case, $Y_{\alpha}(s_{i-1}(\pi_\tau (x)))= \frac{X_{\alpha_1+\cdots 
\alpha_{k-1}}(x) X_{\alpha_{k} + \cdots + \alpha_{i}}(x)}{X_{\alpha_{1}+ \cdots 
\alpha_{i-1}}(x)} = Y_{\alpha}(\pi_{\tau}(s_i x))$.

Case 2: $\alpha= \alpha_ {i-1} + \cdots \alpha_{k-1}$, $i < k$ and 
$w^{-1}(\alpha_i + \cdots \alpha_k) = \tau^{-1}(\alpha_{i-1}+ \cdots 
\alpha_{k-1}) > 0$ and $s_{i-1}\tilde{\tau}=\tau$. 

In this case, $Y_{\alpha}(s_{i-1}(\pi_\tau (x)))= -\frac{X_{\alpha_1+\cdots 
\alpha_{i}}(x) X_{\alpha_{i} + \cdots + \alpha_{k}}(x)}{X_{\alpha_{i}}(x) X_{\alpha_{1}+ 
\cdots \alpha_{k}}(x)} = Y_{\alpha}(\pi_{\tau}(s_i x))$.

Case 3: $\alpha = \alpha_{i-1}$.
$Y_{\alpha}(s_{i-1}(\pi_{\tau}(x)))= \frac{X_{\alpha_{1}+\cdots \alpha_{i}}(x)}
{X_{\alpha_{1}+\cdots \alpha_{i-1}}(x)X_{\alpha_i}(x)}=Y_{\alpha}
(\pi_{\tau}(s_{i-1}(x)))$.

In all other cases, we have:
$Y_{\alpha}(s_{i-1}(\pi_{\tilde{s_{i-1}\tau}}(s_i x)))
=\frac {X_{s_is_1\cdots s_{n}(\alpha)}(x) X_{s_{i}(\beta^{\prime})}(x)}
{X_{s_{i}(s_{1}\cdots s_{n}(\alpha) + \beta^{\prime})}(x)}=Y_{\alpha}(\pi_{\tau}(x))$,  
where $\beta^{\prime}$ is the unique root such that 
$\beta^{\prime} \geq \alpha_1$ and $s_1 \cdots s_n(\alpha)+ \beta^{\prime}$ 
is a root.

This completes the proof.
\end{proof}

With $Y_{\alpha}$'s as in the proof of theorem 5.2, we have 
\begin{corollary}
\beqn
s_1 (Y_{\alpha})= \left\{
\begin{array}{ll} -(1+Y_{\alpha})
 & \mbox{ if } \alpha  \geq \alpha_1. \\
Y_{\alpha} & \mbox{ otherwise } 
\end{array}
\right. .
\eeqn
\end{corollary}
\begin{proof}

Proof follows from the fact that 
\beqn
X_{\alpha}(s_1x)= \left\{\begin{array}{ll} X_{\alpha_1}X_{\alpha}(x)+X_{\alpha_1+\alpha}(x)  & \mbox{ if }  
\alpha = \alpha_2 + \cdots \alpha_i,  ~ 2 \leq i, \\ 
\frac{-X_{\alpha}(x)}{X_{\alpha_1}(x)} & \mbox{ if } \alpha=\alpha_1+ \cdots \alpha_i,  ~ i \geq 2, \\
X_{\alpha}(x) & \mbox{ if } \alpha=\alpha_3+\cdots \alpha_i, ~ i \geq 3
\end{array}
\right. .
\eeqn
\end{proof}
\begin{corollary}
Let $\mathfrak h_n$ be a Cartan subalgebra of $sl_{n+1}(k)$. Let $\chi$ be a 
regular dominant character as in Theorem 5.2. Then, the action of $W$ on the 
GIT quotient 
\[{_{_{T}\backslash \backslash}}({GL_{n+1}(k)/B_{n+1}}))^{ss}(L_{\chi})\simeq 
GL_n(k)/B_n\] is given by the $n$- dimensional representation  $\mathfrak h_n$ of 
$W$.
\end{corollary}
\begin{proof}
Proof follows from theorem 5.2 and corollary 5.3.
\end{proof}

\end{document}